\pdfoutput=1
\RequirePackage{ifpdf}
\ifpdf 
\documentclass[pdftex]{sigma}
\else
\documentclass{sigma}
\fi

\DeclareMathOperator{\Stokesmoduli}{Stokesmoduli}
\DeclareMathOperator{\diag}{diag}
\DeclareMathOperator{\Stm}{Stm}
\DeclareMathOperator{\Mod}{Mod}
\DeclareMathOperator{\Spec}{Spec}
\DeclareMathOperator{\irr}{irr}

\begin{document}

\newcommand{\arXivNumber}{1501.05205}

\allowdisplaybreaks

\renewcommand{\thefootnote}{$\star$}

\renewcommand{\PaperNumber}{036}

\FirstPageHeading

\ShortArticleName{The Stokes Phenomenon and Some Applications}

\ArticleName{The Stokes Phenomenon and Some Applications\footnote{This paper is a~contribution to the Special Issue on
Algebraic Methods in Dynamical Systems.
The full collection is available at
\href{http://www.emis.de/journals/SIGMA/AMDS2014.html}{http://www.emis.de/journals/SIGMA/AMDS2014.html}}}

\Author{Marius VAN DER PUT}
\AuthorNameForHeading{M.~van der Put}
\Address{University of Groningen, Department of Mathematics,\\
P.O.~Box 407, 9700 AK Groningen, The Netherlands}
\Email{\href{mailto:mvdput@math.rug.nl}{mvdput@math.rug.nl}}

\ArticleDates{Received January 21, 2015, in f\/inal form April 21, 2015; Published online May 01, 2015}

\Abstract{Multisummation provides a~transparent description of Stokes matrices which is reviewed here together with some
applications.
Examples of moduli spaces for Stokes matrices are computed and discussed.
A~moduli space for a~third Painlev\'e equation is made explicit.
It is shown that the monodromy identity, relating the topological monodromy and Stokes matrices, is useful for some
quantum dif\/ferential equations and for conf\/luent generalized hypergeometric equations.}


\Classification{14D20; 34M40; 34M55; 53D45}

\renewcommand{\thefootnote}{\arabic{footnote}}
\setcounter{footnote}{0}

\section{Introduction}

Consider a~linear dif\/ferential equation in matrix form $y'+Ay=0$, where the entries of the matrix~$A$ are meromorphic
functions def\/ined in a~neighbourhood of, say, $z=\infty$ in the complex plane.
A~formal or symbolic solution can be lifted to an actual solution in a~sector at $z=\infty$, having the formal solution
as asymptotic behavior.

In 1857, G.G.~Stokes observed, while working in the middle of the night and not long before getting married, the
phenomenon that this lifting depends on the direction of the sector at $z=\infty$ (see~\cite{Ra,St} for more details).

This is the starting point of the long history of the asymptotic theory of singularities of dif\/ferential equations.
The theory of multisummation is the work of many mathematicians such as W.~Balser, B.L.J.~Braaksma, J.~\'Ecalle,
W.B.~Jurkat, D.~Lutz, M.~Loday-Richaud, B.~Malgrange, J.~Martinet, J.-P.~Ramis, Y.~Sibuya (see~\cite{L,LP} for
excellent bibliographies and references and also~\cite{vdP-Si} for some details).

This paper reviews the transparent description of the Stokes phenomenon made possible by multisummation and some
applications of this, namely:
\begin{enumerate}\itemsep=0pt
\item[a)]
moduli spaces of linear dif\/ferential equations,

\item[b)] quantum dif\/ferential equations and conf\/luent generalized hypergeometric equations,

\item[c)] isomonodromy, the Painlev\'e equations and Okamoto--Painlev\'e spaces.
\end{enumerate}

Moreover the paper presents new ideas on moduli spaces for the Stokes phenomenon and the use of the monodromy identity.

The f\/irst section is written for the convenience of the reader.
It is a~review of a~large part of~\cite{vdP-Si} and has as new part Section~\ref{section1.5}.
The aim is to describe clearly and concisely the theory and results while bypassing technical details and proofs.

One considers a~singular matrix dif\/ferential equation $y'+Ay=0$ at $z=\infty$ and recalls its formal classif\/ication, the
def\/inition of the Stokes matrices and the analytic classif\/ication.
The theory is illustrated in Section~\ref{section1.5} by the conf\/luent generalized hypergeometric equation~$_pD_q$.
We rediscover results from~\cite{DM,Mi} as application of the monodromy identity.
In the second section moduli spaces for Stokes maps are discussed.
In the new part Section~\ref{section2.2} and especially Proposition~\ref{proposition2}, it is shown that totally ramif\/ied equations have very interesting  Stokes matrices.
Quantum dif\/ferential equations coming from Fano varieties are studied in Section~\ref{section3}.
Again the computability of Stokes matrices is the theme.
Certain moduli spaces, namely Okamoto--Painlev\'e spaces, corresponding to Painlev\'e equations are discussed in Section~\ref{section4},
including an explicit calculation of a~monodromy space for ${\rm P_{III}(D_7)}$.

\section{Formal and analytic classif\/ication, Stokes maps}\label{section1}

\subsection{Terminology and notation}\label{section1.1}

Let~$k$ be a~{\it differential field}, i.e., a~f\/ield with a~map $f\mapsto f'$ (called a~derivation) satisfying $(f+g)'=f'+g'$ and $(fg)'=f'g+fg'$.
The f\/ield of constants of~$k$ is $\{f\in k\,|\, f'=0\}$.
In this paper we suppose that the f\/ield of constants of~$k$ is $\mathbb C$ and that $k\neq \mathbb C$.

The most important dif\/ferential f\/ields that we will meet are $\mathbb{C}(z)$, $\widehat{K}:=\mathbb{C}((z^{-1}))$ and
$K:=\mathbb{C}(\{z^{-1}\})$, i.e., the f\/ield of the rational functions in~$z$, the formal Laurent series in $z^{-1}$ and
the f\/ield of the convergent Laurent series in $z^{-1}$.
The latter is the f\/ield of germs of meromorphic functions at $z=\infty$.
In all cases the dif\/ferentiation is $f\mapsto f'=\frac{df}{dz}$ (sometimes replaced by $f\mapsto
\delta(f):=z\frac{df}{dz} $ in order to make the formulas nicer).

A matrix dif\/ferential equation $y'+Ay=0$ with~$A$ a~$d\times d$-matrix with coordinates in~$k$ gives rise to the
operator $\partial:=\frac{d}{dz}+A$ of $k^d$, where $\frac{d}{dz}$ acts coordinatewise on $k^d$ and~$A$ is the matrix of
a~$k$-linear map $k^d\rightarrow k^d$.
Write now $M=(M,\partial)$ for $k^d$ and the operator $\partial$.
Then this object is a~dif\/ferential module.

Indeed, a~{\it differential module} over~$k$ is a~f\/inite dimensional~$k$-vector space~$M$ equipped with an additive map
$\partial\colon M\rightarrow M$ satisfying $\partial (fm)=f'm+f\partial (m)$ for any $f\in k$ and $m\in M$.
If one f\/ixes a~basis of~$M$ over~$k$, then~$M$ is identif\/ied with $k^d$ (with $d=\dim M$) and the ope\-ra\-tor~$\partial$ is
identif\/ied with $\frac{d}{dz}+A$.
Here~$A$ is the matrix of $\partial$ with respect to the given basis of~$M$.

Thus a~dif\/ferential module is ``a matrix dif\/ferential equation where the basis is forgotten'' and a~matrix dif\/ferential  
equation is the same as a~dif\/ferential module with a~given basis.

We will use dif\/ferential operators, i.e., elements of the skew polynomial ring $k[\partial]$ (where $\partial$ stands
for $\frac{d}{dz}$) def\/ined by the rule $\partial f=f\partial +f'$.
Instead of $\partial$ we sometimes use $\delta:=z\partial$.
Then $\delta f=f\delta +\delta (f)$ (in particular $\delta z=z\delta +z$).

Let~$M$ be a~dif\/ferential module.
The ring $k[\partial]$ acts from the left on~$M$.
For any element $e\in M$, there is a~monic operator $L\in k[\partial]$ of smallest degree such that $Le=0$.
The element~$e$ is called {\it cyclic} if~$L$ has degree $d=\dim M$.
Cyclic elements~$e$ exist \cite[\S\S~2.10 and~2.11]{vdP-Si}.
The corresponding operator~$L$ is identif\/ied with a~scalar dif\/ferential equation and~$L$ determines~$M$.

In practice one switches between dif\/ferential modules, matrix dif\/ferential equations, scalar dif\/ferential equations and
dif\/ferential operators.

\subsection[Classif\/ication of dif\/ferential modules over $\widehat{K}:=\mathbb{C}((z^{-1}))$]{Classif\/ication
of dif\/ferential modules over $\boldsymbol{\widehat{K}:=\mathbb{C}((z^{-1}))}$}\label{section1.2}

The classif\/ication of a~matrix dif\/ferential equation $z\frac{d}{dz}+A$ over $\widehat{K}$ (note that we prefer here~$z\frac{d}{dz}$) goes back to G.~Birkhof\/f and H.L.~Turritin.
This classif\/ication is somewhat similar to the Jordan normal of a~matrix.
However it is more subtle since $z\frac{d}{dz}+A$ is linear over~$\mathbb C$ and is not linear over~$\widehat{K}$.

We prefer to work ``basis free'' with a~dif\/ferential module~$M$ and classify~$M$ by its {\it solution space~$V$} with
additional data forming a~tuple $(V,\{V_q\},\gamma)$.

If $\dim M=d$, then we want the solution space, i.e., the elements~$w$ with $\delta w=0$ to be a~$\mathbb C$-vector
space of dimension~$d$.
Now we write $\delta:M\rightarrow M$ instead of $\partial$ because $\frac{d}{dz}$ is replaced by $z\frac{d}{dz}$.

In general $\{m\in M \,|\, \delta(m)=0\}$ is a~vector space over $\mathbb C$ with dimension $<d$.
Thus we enlar\-ge~$\widehat{K}$ to a~suitable dif\/ferential ring~$U$ and consider $\{w\in U\otimes_{\widehat{K}} M \,|\, \delta w=0\}$.

This dif\/ferential ring~$U$ (the universal Picard--Vessiot ring for $\widehat K$) is built as follows.
We need a~linear space of ``{\it eigenvalues}'' $\mathcal{Q}:=\cup_{m\geq 1}z^{1/m}\mathbb{C}[z^{1/m}]$ and {\it
symbols} $z^\lambda$ with $\lambda \in \mathbb{C}$, $\log z$, $e(q)$ with $q\in \mathcal{Q}$.
The relations are $z^{a+b}=z^a\cdot z^b$, $z^1=z\in \widehat{K}$, $e(q_1+q_2)=e(q_1)\cdot e(q_2)$, $e(0)=1$.
And we def\/ine their derivatives by the formulas $(z^a)'=az^a$, $\log (z)'=1$, $e(q)'=qe(q)$ (we note that $'$ stands for
$z\frac{d}{dz}$ and that the interpretation of $e(q)$ is $e^{\int q\frac{dz}{z}}$).

This {\it universal Picard--Vessiot ring} is $U:=\widehat{K}[\{z^\lambda\}_{\lambda \in \mathbb{C}}, \log(z),
\{e(q)\}_{q\in \mathcal{Q}}]$.
This ring is a~direct sum $U=\oplus_{q\in \mathcal{Q}}U_q$ and $U_q:=e(q) \widehat{K}[\{z^\lambda\}_{\lambda \in
\mathbb{C}}, \log(z)]$.

The Galois group of the algebraic closure $\cup_{m\geq 1}\mathbb{C}((z^{-1/m}))$ of $\widehat{K}$ is $\cong
\widehat{\mathbb{Z}}$ and is topologically generated by the element~$\gamma$ given by $\gamma z^\lambda =e^{2\pi i
\lambda}z^\lambda$ for all $\lambda \in \mathbb{Q}$.
The algebraic closure of~$\widehat{K}$ lies in~$U$.
One extends~$\gamma$ to a~dif\/ferential automorphism of~$U$ by the following formulas (corresponding to the
interpretation of the symbols) $\gamma z^a=e^{2\pi i a}z^a$ for all $a\in \mathbb{C}$, $\gamma \log(z)= 2\pi i
+\log(z)$, $\gamma e(q)=e(\gamma q)$.

For every dif\/ferential module~$M$ over $\widehat{K}$, its solution space, def\/ined as $V:=\ker (\delta,U\otimes M)$ has
``all solutions'' in the sense that $\dim_\mathbb{C}V=\dim_{\widehat{K}}M$ and the canonical map $U\otimes_\mathbb{C}V
\rightarrow U\otimes_{\widehat{K}}M$ is an isomorphism.

Put $V_q:=\ker (\delta,U_q\otimes M)$.
Then $V=\oplus_qV_q$ is a~decomposition of the solution space.
Further the action of~$\gamma$ on~$U$ induces a~$\gamma \in {\rm GL}(V)$ such that $\gamma V_q=V_{\gamma q}$ for all~$q$.

\begin{theorem}[\protect{formal classif\/ication \cite[Proposition~3.30]{vdP-Si}}]  \label{theorem1}
The functor $M\mapsto (V,\{V_q\},\gamma)$ is an
equivalence of the Tannakian categories of the differential modules over $\widehat{K}$ and the category of the tuples
$(V,\{V_q\},\gamma)$.
\end{theorem}

The category of the tuples $(V,\{V_q\},\gamma)$ is denoted by ${\rm Gr}_1$ in~\cite[p.~76]{vdP-Si}.
The ``Tannakian'' property of the functor of the theorem means that it commutes with all constructions of linear
algebra, including tensor products, applied to modules.

Suppose that~$M$ induces the tuple $(V,\{V_q\},\gamma)$.
Then~$q$ is called an {\it eigenvalue of~$M$} if $V_q\neq 0$.
{\it The Katz invariant of~$M$} is the maximum of the degrees in~$z$ of the eigenvalues~$q$ of~$M$.
Clearly, the Katz invariant is a~non negative rational number which measures the singularity.
In particular, the Katz invariant is zero if and only if the module~$M$ is regular singular.

Further the map $\gamma \in {\rm GL}(V)$ is called the {\it formal monodromy of~$M$}.
Section~\ref{section1.3} illustrates the computation of the tuple $(V,\{V_q\},\gamma)$.

\subsection{The conf\/luent generalized hypergeometric equation}\label{section1.3}

Consider the conf\/luent generalized hypergeometric equation in operator form
\begin{gather*}
_pD_q=(-1)^{q-p}z\prod\limits_{j=1}^p(\delta +\mu_j)- \prod\limits_{j=1}^q(\delta +\nu_j-1)
\qquad
\text{with}
\qquad
\delta=z\frac{d}{dz}.
\end{gather*}
We assume that $1\leq p<q$ and that the complex parameters $\mu_j$, $\nu_j$ are such that $\mu_1,\dots,\mu_p$ are
distinct modulo~$\mathbb{Z}$.
We regard the operator~$_pD_q$ as element of~$\widehat{K}[\delta]$.

Its Newton polygon has two slopes $0$ and $\frac{1}{q-p}$.
The operator has then for each ordering of the slopes a~unique factorization (compare~\cite[Theorem 3.8]{vdP-Si}).
The two decompositions are
\begin{gather*}
_pD_q=-\big(\delta^p+a_{p-1}\delta^{p-1}+\dots +a_0\big) \big(\delta^{q-p}+b_{q-p-1}\delta^{q-p-1}+\dots +b_1\delta-(-1)^{q-p}z+b_0\big)
\end{gather*}
with all $a_j, b_j\in \mathbb{C}[[z^{-1}]]$.
The two factors almost commute and there is a~factorization in the opposite ordering of the slopes
\begin{gather*}
_pD_q=-\big(\delta^{q-p}+b^*_{q-p-1}\delta^{q-p-1}+\dots + b^*_1\delta -(-1)^{q-p}z+b^*_0\big)
\big(\delta^p+a^*_{p-1}\delta^{p-1}+\dots +a^*_0\big)
\end{gather*}
with all $a^*_j, b^*_j\in \mathbb{C}[[z^{-1}]]$.

From these factorizations one can read of\/f the solution space.
The term $(\delta^p+a^*_{p-1}\delta^{p-1}+\dots +a^*_1\delta +a^*_0)$ is {\it equivalent} to $\prod\limits_{j=1}^p(\delta
+\mu_j)$ (`equivalent' means that the dif\/ferential modules over $\widehat{K}$ def\/ined by the two operators are
isomorphic).
This yields solutions $f_j:=z^{-\mu_j}$ for $j=1,\dots,p$.
They form a~basis of the $\mathbb C$-vector space $V_0$ with {\it eigenvalue $q_0=0$}.
The action of the formal monodromy~$\gamma$ on $V_0$ is $\gamma (f_j)=e^{-2\pi i\mu_j}f_j$.

The term $(\delta^{q-p}+b_{q-p-1}\delta^{q-p-1}+\dots +b_1\delta -(-1)^{q-p}z+b_0)$ is equivalent to $\delta^{\sigma}-z$
(up to a~sign) with $\sigma:=q-p$.
This operator factors over $\widehat{K}(z^{1/\sigma})$ and one f\/inds the {\it eigenvalues} $q_1=z^{1/\sigma}$,
$q_2=\zeta z^{1/\sigma}$, $\dots$, $q_\sigma=\zeta^{\sigma -1}z^{1/\sigma}$ with $\zeta =e^{2\pi i/\sigma}$.

{\it Now we can describe the solution space~$V$ of $_pD_q$}:
\begin{gather*}
V=V_0\oplus V_{q_1}\oplus \dots \oplus V_{q_\sigma},
\end{gather*}
where $V_0$ has a~basis $f_1,\dots,f_p$ with $\gamma f_j=e^{-2\pi i \mu_j}f_j$.
Choose a~basis $e_1$ of the 1-dimensional space $V_{q_1}$.
Put $e_2:=\gamma e_1$, $e_3:=\gamma e_2$, $\dots$, $e_\sigma:=\gamma e_{\sigma -1}$.
Then $V_{q_j}=\mathbb{C}e_j$ for $j=1,\dots,\sigma$.
Finally $\gamma e_\sigma = e^{2\pi i \lambda} e_1$ with $\lambda:=\frac{1}{2}(\sigma +1)+\sum\limits_{j=1}^p\mu_j
-\sum\limits_{j=1}^q\nu_j$.
This follows from a~computation of~$\gamma$ on $f_1\wedge \dots \wedge f_p\wedge e_1\wedge \dots \wedge e_\sigma$ and
the determinant $\delta +\sum\limits_{j=1}^q(v_j-1)$ of $_pD_q$.

The above coincides with the formula of the formal monodromy in~\cite[p.~373]{Mi}.
In~\cite{DM,Mi} an explicit basis of formal or symbolic solutions of $_pD_q$ is constructed and the
computation of the formal monodromy and, later on, of the Stokes matrices is with respect to this basis.
Our basis $f_1,\dots,f_p,e_1,\dots,e_\sigma$ is not unique.
More precisely, the above tuple has a~non trivial automorphism group $G\cong (\mathbb{C}^*)^{p+1}$.
The elements of~$G$ are given by $f_j\mapsto \alpha_jf_j$ for $j=1,\dots,p$ and $e_j\mapsto \alpha_{p+1}e_j$ for
$j=1,\dots,\sigma $ (and $(\alpha_1,\dots,\alpha_{p+1})\in G$).
We will see that the group $G/\mathbb{C}^*$ acts non trivially on the entries of the Stokes maps (see also
Proposition~\ref{proposition2}).

\subsection{Stokes maps and the analytic classif\/ication}\label{section1.4}

Let~$M$ be a~dif\/ferential module over~$K$.
Then $\widehat{K}\otimes M$ is a~dif\/ferential module over $\widehat{K}$ and induces a~tuple $(V,\{V_q\},\gamma)$.
For two eigenvalues $q$, $\tilde{q}$ of $\widehat{K}\otimes M$ one considers special directions $e^{2\pi i d}$, $d\in \mathbb{R}$,
called {\it singular} for the dif\/ference $q-\tilde{q}$.
Those are the~$d$ such that $e^{\int (q-\tilde{q})\frac{dz}{z}}$ (this is the solution of $y'=(q-\tilde{q})y$) has
maximal descent to zero for $z:=re^{2\pi i d}$ and $r>0$, $r\rightarrow 0$.

The tuple $(V,\{V_q\},\gamma)$ is the formal classif\/ication of~$M$, i.e., the classif\/ication of $\widehat{K}\otimes_KM$.
Now we consider the classif\/ication of~$M$ itself.
In the sequel, we use {\it multisummation} as a~{\it black box}.
It is a~technical extension of the classical Borel summation of certain divergent power series.
We refer to~\cite{B1,B2,vdP-Si} for details.

However, we will not need {\it the precise formulation of the multisum in a~certain direction~$d$}.
This is rather technical \cite[Def\/inition 7.46]{vdP-Si}.

{\it The sloppy definition and result is as follows}.
For a~direction~$d$ which is not singular, one considers a~certain sector $S(d)$ at $z=\infty$ around~$d$ (or more
precisely, a~sequence of nested sectors).
There is a~{\it unique} $\mathbb{C}$-linear map, the multisum
\begin{gather*}
\operatorname{sum}_d\colon \  V\rightarrow \text{the space of the solutions of $M$ on the sector $S(d)$}.
\end{gather*}
It has the property that for any $v\in V$, the ${\rm sum}_d(v)$ has asymptotic expansion~$v$ on the prescribed sector
(or more precisely, a~sequence of asymptotic properties on the nested sequence of sectors).

We note that the uniqueness of $\operatorname{sum}_d$ is the important issue of multisummation.
This is a~much stronger result than the more classical ``Main asymptotic existence theorem'' \cite[Theorem~7.10]{vdP-Si}.

For a~singular direction~$d$ one takes real numbers $d^-<d<d^+$ close to~$d$ and def\/ines the Stokes map $\operatorname{St}_d\in {\rm
GL}(V)$ by $\operatorname{sum}_{d^+}=\operatorname{sum}_{d^-}\circ \operatorname{St}_d$.
To~$M$ we associate the tuple $(V,\{V_q\},\gamma, \{\operatorname{St}_d\}_{d\in \mathbb{R}})$.

One can show \cite[Theorem~8.13 and Remark 8.14]{vdP-Si} that this tuple has the additional pro\-per\-ties:
\begin{enumerate}\itemsep=0pt
\item[$(*)$] $\operatorname{St}_d$ has the form $\operatorname{id}+\sum\limits_{d \; \text{singular for} \; q-\tilde{q}} \operatorname{Hom}(V_q,V_{\tilde{q}})$.
Here, every $\ell\in \operatorname{Hom}(V_q,V_{\tilde{q}})$ is seen as an element of $\operatorname{End}(V)$ by the sequence of maps
$V\stackrel{\rm projection}{\rightarrow}V_q\stackrel{\ell}{\rightarrow}V_{\tilde{q}} \stackrel{\rm
inclusion}{\rightarrow}V$.

\item[$(**)$] $\gamma^{-1}\operatorname{St}_d\gamma =\operatorname{St}_{d+1}$.
\end{enumerate}

One considers the category of the tuples with the properties $(*)$ and $(**)$.
This category of f\/inite dimensional complex vector spaces with this additional ``linear algebra structure'' is denoted
by ${\rm Gr}_2$ in~\cite[\S~9.2]{vdP-Si}.
It is a~Tannakian category.

\begin{theorem}[\protect{the analytic classification \cite[Theorem~9.11]{vdP-Si}}]\label{theorem2}
The functor $M\!\mapsto\! (V,\{V_q\},\gamma, \{\operatorname{St}_d\})$ is an equivalence of the Tannakian categories of the differential
modules over~$K$ and the above category of the tuples $(V,\{V_q\},\gamma, \{\operatorname{St}_d\})$ satisfying $(*)$ and $(**)$.
\end{theorem}

The above theorem is, in contrast with the formal classif\/ication, {\it a~deep and final result in the asymptotic theory
of linear differential equations}.

The {\it irregularity of Malgrange}, $\irr(M)$ of the dif\/ferential module~$M$ is def\/ined by $\irr(M)=\sum\limits_{q\neq
\tilde{q}}\deg_{z}(q-\tilde{q})\cdot \dim V_q\cdot \dim V_{\tilde{q}}$.
One observes that the dimension of the space of all possibilities for Stokes maps with a~f\/ixed formal tuple is equal to
$\irr(M)$ (see also Section~\ref{section2}).
A~useful result, obtained by the above description of the Stokes maps, is the following.

\begin{proposition}[\protect{the monodromy identity \cite[Proposition~8.12]{vdP-Si}}]\label{proposition1}
The topological monodromy of~$M$ is conjugated to $\gamma \operatorname{St}_{d_s}\cdots \operatorname{St}_{d_1}\in {\rm GL}(V)$,
where $0\leq d_1<\dots <d_s<1$ are the singular directions of~$M$.
\end{proposition}

One cannot claim that the topological monodromy is equal to this product since one has to identify~$V$ with the local
solution space at a~point near $z=\infty$ and that can be done in many ways.

\subsection[Stokes matrices for $_pD_q$]{Stokes matrices for $\boldsymbol{_pD_q}$}\label{section1.5}

The {\it main observation} is the following.
The monodromy identity yields {\it complete formulas for the Stokes matrices} of $_pD_q$ if we are allowed to choose
a~suitable basis of~$V$.
In this way we rediscover the formulas for the Stokes matrices of~\cite{DM,Mi}.
More precisely, by choosing multiples of the basis $f_1,\dots,f_p$, we can normalize~$p$ entries of the Stokes maps to be~1.
The others are then determined by the monodromy identity.
This works, according to Proposition~\ref{proposition2}, under the assumption that $\mu_1,\dots,\mu_p$ are distinct modulo $\mathbb{Z}$
and that the equation is irreducible.

We note that {\it the differential Galois group of $_pD_q$} is in fact the dif\/ferential Galois group of $_pD_q$ as
equation over the f\/ield of convergent Laurent series~$K$.
This group does not depend on the choice of multiples of $f_1,\dots,f_p $.
For the rather involved computation of this dif\/ferential Galois group, see~\cite{DM, Mi}, the knowledge of the formal
classif\/ication and the characteristic polynomial of the topological monodromy at $z=0$ (or equivalently at $z=\infty$)
suf\/f\/ice.
We illustrate this by the easy example $_1D_3=z(\delta +\mu)-\prod\limits_{j=1}^3(\delta +\nu_j-1)$.

The formal solution space is the direct sum of three 1-dimensional spaces $V_0\oplus V_{z^{1/2}}\oplus V_{-z^{1/2}}$
with basis $f_1$, $e_1$, $e_2$.
There is only one singular direction in the interval $[0,1)$, namely $d=0$.
The topological monodromy at $z=0$ is that of the operator $\prod\limits_{j=1}^3(\delta +\nu_j-1)$ and has eigenvalues
$e^{-2\pi i \nu_j}$, $j=1,2,3$.

The formal monodromy is $\gamma (f_1)=e^{- 2\pi i \mu}f_1$, $\gamma e_1=e_2$ and $\gamma (e_2)= e^{-2\pi i \lambda}
e_1$.
The product of the formal monodromy and the unique Stokes matrix in the interval $[0,1)$ (its form is described in Theorem~\ref{theorem2}) is
\begin{gather*}
\left(
\begin{matrix}
 d_1 & 0 &0
\\
0&0 &d_2
\\
0 &1 &0
\end{matrix}
\right) \left(
\begin{matrix}
 1 &x_{1,0} & 0
\\
0&1 &0
\\
x_{0,2} &x_{1,2} & 1
\end{matrix}
\right)= \left(
\begin{matrix} d_1 &d_1x_{1,0} &0
\\
d_2x_{0,2} &d_2x_{1,2} & d_2
\\
0 &1 & 0
\end{matrix}
\right),
\end{gather*}
where $d_1=e^{-2\pi i \mu}$, $d_2=e^{-2\pi i \lambda}$, $\lambda =\frac{3}{2}+\mu -\sum \nu_j$.
Its characteristic polynomial $T^3-(d_1+d_2x_{1,2})T^2-(d_1d_2x_{1,0}x_{0,2}-d_2)T+d_1d_2$ coincides with
$\prod\limits_{j=1}^3(T-e^{2\pi i v_j})$.
Choose $x_{1,0}=1$.
Then all Stokes matrices are determined.

The exceptional case $x_{1,0}x_{0,2}=0$ cannot be handled in this way.
Using J.-P.~Ramis result \cite[Theorem~8.10]{vdP-Si} and~\cite[Th\'eor\`eme~4.10]{CR} that the dif\/ferential Galois
group is generated as algebraic group by the formal monodromy, the exponential torus and the Stokes matrices, one
concludes that $x_{1,0}=0$ or $x_{0,2}=0$ implies that the equation is reducible (compare the proof of Proposition~\ref{proposition2}).
See~\cite{DM} for a~complete description of all cases and all dif\/ferential Galois groups.
We remark that the monodromy identity for~$_pD_q$ is explicitly present in~\cite{DM}.

\section{Moduli spaces for the Stokes data}\label{section2}

For given formal data $F:=(V,\{V_q\},\gamma)$ at $z=\infty$, there exists a~unique dif\/ferential module~$N$ over
$K=\mathbb{C}(\{z^{-1}\})$ with these formal data and with trivial Stokes matrices.
One considers dif\/ferential modules~$M$ which have formal classif\/ication~$F$.
The set of isomorphism classes of these modules {\it does not have a~good algebraic structure} since~$F$ has, in
general, automorphisms.

D.G.~Babbitt and V.S.~Varadarajan~\cite{BV} consider instead pairs $(M,\phi)$ of a~dif\/ferential mo\-du\-le~$M$ over~$K$
and an isomorphism $\phi:\widehat{K}\otimes M\rightarrow \widehat{K}\otimes N$.
Two pairs $(M_j,\phi_j)$, $j=1,2$, are equivalent if there exists an isomorphism $\alpha:M_1\rightarrow M_2$ such that $\phi_2\circ\alpha=\phi_1$.
The set $\Stokesmoduli(F)$ of equivalence classes of pairs $(M,\phi)$ has been given a~natural structure of complex algebraic variety.
Babbitt and Varadarajan prove that $\Stokesmoduli(F)$ is isomorphic to the af\/f\/ine space $\mathbb{A}_\mathbb{C}^m$, where
$m=\sum\limits_{i\neq j}\dim V_{q_i}\cdot \dim V_{q_j} \cdot \deg_{z}(q_i-q_j)$.
The $\{q_i\}$ are the eigenvalues of~$M$ and~$m$ is the irregularity of~$F$, in the terminology of B.~Malgrange.

This result of is in complete agreement with the above description of the Stokes matrices~$\{\operatorname{St}_d\}$ (def\/ined~by multisummation).
Thus $\Stokesmoduli(F)$ is the moduli space for the possible Stokes matrices for f\/ixed formal data $F=(V,\{V_q\},\gamma)$.

However, there is, {\it in general}, no universal family of dif\/ferential modules parametrized by $\Stokesmoduli(F)\cong
\Spec(\mathbb{C}[x_1,\dots,x_{m}])$.
In other words, $\Stokesmoduli(F)$ {\it is, in general, not a~fine moduli space for the above family of differential modules}.

Indeed, suppose that such a~family $\{M_\xi \,|\, \xi\in \Stokesmoduli(F)\}$ of dif\/ferential modules over $K=\mathbb{C}(\{z^{-1}\})$ exists.
This family is represented by a~matrix dif\/ferential operator $z\frac{d}{dz} +A$ in the variable~$z$ and with entries in,
say, $K(x_1,\dots,x_{m})$.
The monodromy identity shows that the eigenvalues of the topological monodromy are algebraic over this f\/ield.
A~logarithm of the topological monodromy is computable from $z\frac{d}{dz} +A$ and has again entries in $K(x_1,\dots,x_{m})$.
This is, {\it in general}, not possible.
See Section~\ref{section2.1} for a~concrete case.

In order to produce a~f\/ine moduli space one replaces the dif\/ferential module~$M$ over $K=\mathbb{C}(\{z^{-1}\})$~by
a~tuple $(\mathcal{M},\nabla,\phi)$.
Here $(\mathcal{M},\nabla)$ is a~connection on the projective line $\mathbb{P}^1$ over $\mathbb{C}$ which has two
singular points $0$ and~$\infty$.
The point $z=0$ is supposed to be regular singular.
Further~$\phi$ is an isomorphism of the formal completion of the connection at $z=\infty$ (i.e., $\mathcal{M}_\infty
\otimes \mathbb{C}[[z^{-1}]]$) with a~prescribed object $\nabla:N_0\rightarrow N_0\otimes z^k\mathbb{C}[[z^{-1}]]$ (for
suitable~$k$, see below).

This prescribed object is the following.
Let~$N$ denote the dif\/ferential module over $\mathbb{C}((z^{-1}))$ corresponding to the given $F=(V,\{V_q\},\gamma)$.
Then $N_0\subset N$ is the $\mathbb{C}[[z^{-1}]]$-submodule of~$N$, generated by a~basis of~$N$ over
$\mathbb{C}((z^{-1}))$ (i.e., $N_0$ is a~lattice in~$N$).
There are many choices for $N_0$ and in principle the choice is not important.
However we choose a~standard lattice $N_0$ (see~\cite[\S\S~12.4 and~12.5]{vdP-Si}) in order to make explicit computation easier.
The given dif\/ferential operator $\nabla:N\rightarrow N$ maps $N_0$ into $z^kN_0\subset N$ for a~certain integer $k\geq0$.
The choice of~$k$ is irrelevant for the construction.

According to a~theorem of Birkhof\/f \cite[Lemma~12.1]{vdP-Si}, this ``spreading out of~$M$'' exists.
We will make the assumption that $\mathcal{M}$ is a~{\it free vector bundle} on $\mathbb{P}^1$.
The above description leads to a~f\/ine moduli space $\Mod(F)$.

Let $\Stm: \Mod(F)\rightarrow \Stokesmoduli(F)$, denote the map which associates to a~tuple $(\mathcal{M},\nabla, \phi)$,
belonging to $\Mod(F)$, its set of Stokes matrices.
Known results are:

\begin{theorem}[\protect{\cite[\S\S~12.11, 12.17, 12.19, 12.20]{vdP-Si}}] \label{theorem3}\quad
\begin{enumerate}  \itemsep=0pt
\item[$(a)$]
$\Mod(F)$ is isomorphic to the affine space $\mathbb{A}_\mathbb{C}^m$.
\item[$(b)$] $\Stm$ is analytic and has an open dense image.

\item[$(c)$] The generic fibre of $\Stm$ is a~discrete infinite set and can be interpreted as a~set of logarithms of the topological monodromy.
\end{enumerate}
\end{theorem}

{\it Comments}.
The proof of (a) is complicated and the result itself is somewhat amazing.

(b) follows from the observations: If the topological monodromy of~$M$ is semi-simple, then a~tuple
$(\nabla,\mathcal{M},\phi)$ with free $\mathcal{M}$ exists.
Moreover semi-simplicity is an open property.

(c) follows from the construction of ``spreading out''.
One needs a~logarithm of the topological monodromy in order to construct the connection $(\mathcal{M},\nabla)$ on
$\mathbb{P}^1$ from the dif\/ferential module over $K=\mathbb{C}(\{z^{-1}\})$.

A precise description of the f\/ibres seems rather dif\/f\/icult.
Moreover, a~better moduli space, replacing $\Mod(F)$, which does not require the vector bundle $\mathcal{M}$ to be free,
should be constructed.

\subsection{Example: unramif\/ied cases}\label{section2.1}

$F$ is def\/ined by $V=V_{\lambda_1z}\oplus \dots \oplus V_{\lambda_nz}$, where each $V_{\lambda_jz}$ has dimension one
and the $\lambda_1,\dots,\lambda_n\in \mathbb{C}$ are distinct.
Further $\gamma $ is the identity.
Then $N_0$ can be given by the dif\/ferential operator $\delta +z \cdot \diag(\lambda_1,\dots,\lambda_n)$.
Clearly $\Stokesmoduli(F)\cong \mathbb{A}_\mathcal{C}^{n(n-1)}$.

The universal family is $\delta +z \cdot \diag(\lambda_1,\dots,\lambda_n)+ (T_{i,j})$, where for notational convenience
$T_{i,i}=0$ and the $\{T_{i,j}\}$ with $i\neq j$ are $ n^2-n$ independent variables \cite[Theorem~12.4]{vdP-Si}.
Thus $\Mod(F)$ is indeed isomorphic to $\mathbb{A}_\mathbb{C}^{n(n-1)}$.
Further one observes that, in general, the matrix $L:=(T_{i,j})$ has the property that $e^{2\pi iL}$ is the topological
monodromy at $z=0$ (or equivalently at $z=\infty$).
Now, by the monodromy identity (Proposition~\ref{proposition1}), the entries of $e^{2\pi i L}$ are (up to conjugation) rational in the
$n(n-1)$-variables of $\Stokesmoduli(F)$.
This shows that there is no universal family above $\Stokesmoduli(F)$.

In the case $n=2$, the map $\Mod(F)\rightarrow \Stokesmoduli(F)$ can be made explicit and is shown to be surjective.
For $n>2$, the above map is ``highly transcendental'' and we do not know whether it is surjective.
The problem is the choice of a~free vector bundle $\mathcal{M}$ in the def\/inition of $\Mod(F)$.

The problem of explicit computation of the Stokes matrices, i.e., making $\Stm$ explicit in this special case, has been
studied over a~long period and by many people G.D.~Birkhof\/f, H.~Turrittin, W.~Balser, W.B.~Jurkat, D.A.~Lutz, K.~Okubo,
B.~Dubrovin, D.~Guzzetti et al.\
(see the introduction of~\cite{BJL1} and also~\cite{BJL2, Du98, Guz}).
Using Laplace integrals the solutions of the above equation can be expressed in solutions of a~`transformed equation'
with only regular singularities.
The ordinary monodromy of the transformed equation produces answers for the Stokes matrices of the original equation.
This is used by Dubrovin (see~\cite[Lemma~5.4, p.~97]{Dub99}), D.~Guzzetti, H.~Iritani, S.~Tanabe, K.~Ueda et al.
in computations of the Stokes matrix for quantum dif\/ferential equations, see~\cite{Dub96, Du98, Dub99, Guz, Iri, Tan, T-U, Ue1,Ue2}.

{\it Now we come to a~surprising new case}.

\subsection{Example: Totally ramif\/ied cases}\label{section2.2}

The formal data~$F$ is essentially $V=V_{z^{1/n}}\oplus V_{\zeta z^{1/n}}\oplus \dots \oplus V_{\zeta^{n-1}z^{1/n}}$,
where $\zeta:=e^{2\pi i/n}$, each $V_{\zeta^jz^{1/n}}$ has dimension 1 and~$\gamma$ satisf\/ies $\gamma^n=1$.
The irregularity $\sum\limits_{i\neq j} 1\cdot 1\cdot \deg (\zeta^iz^{1/n}-\zeta^jz^{1/n})$ is equal to $n-1$ and this
is small compared to the unramif\/ied case with irregularity $n(n-1)$.
This is responsible for special features of these important examples.

{\it For notational convenience we will consider the case $n=3$}.
\\
The lattice $N_0$ with formal data~$F$ and trivial Stokes matrices can be represented by the dif\/ferential operator
$\delta +\left(
\begin{matrix} -\frac{1}{3} &0 &z
\\
1 &0 &0
\\
0 &1 & \frac{1}{3}
\end{matrix}
\right)$.
A~computation (following~\cite[\S~12.5]{vdP-Si}) shows that $\Mod(F)$ is represented by the universal family $\delta+\left(
\begin{matrix}
a_1&0 &z
\\
1 &a_2 &0
\\
0 &1 & a_3
\end{matrix}
\right) $ with $a_1,a_2,a_3\in\mathbb{C}$ with $a_1+a_2+a_3=0$.

There are 6 singular directions, corresponding to the dif\/ferences of generalized eigenvalues
$\zeta^iz^{1/3}-\zeta^jz^{1/3}$ for $i\neq j$.
The corresponding Stokes matrix has one element of\/f the diagonal, called $x_{j,i}$.
Two singular directions are in $[0,1)$, namely $\frac{1}{4}$, $\frac{3}{4}$, and the others are obtained by shifts over $1$ and~$2$.
The topological monodromy at $z=0$ is conjugated (by the monodromy identity) to $\gamma \operatorname{St}_{3/4}\operatorname{St}_{1/4}$ which reads
\begin{gather*}
\left(
\begin{matrix} 0&0 &1
\\
1 &0 &0
\\
0&1 &0
\end{matrix}
\right) \left(
\begin{matrix}1 &0 &0
\\
0&1 &0
\\
0 &x_{2,1} &1
\end{matrix}
\right) \left(
\begin{matrix}1 &x_{0,1} &0
\\
0 &1 &0
\\
0 &0 & 1
\end{matrix}
\right).
\end{gather*}
The characteristic polynomial of this matrix is $\lambda^3-x_{0,1}\lambda^2-x_{2,1}\lambda -1$.
The topological monodromy is given by the operator $\delta + \left(
\begin{matrix} a_1&0 &0
\\
1 &a_2 &0
\\
0 &1 & a_3
\end{matrix}
\right)$,
which has eigenvalues $a_1$, $a_2$, $a_3$.

The monodromy has eigenvalues $e^{2\pi i a_j}$ for $j=1,2,3$ and its characteristic polynomial is
\begin{gather*}
\big(\lambda - e^{2\pi i a_1}\big)\big(\lambda -e^{2\pi i a_2}\big)\big(\lambda -e^{2\pi i a_3}\big).
\end{gather*}
Hence $x_{0,1}=\sum\limits_{j=1}^3e^{2\pi i a_j}$ and $x_{2,1}=\sum\limits_{i<j} e^{2\pi i a_i} e^{2\pi i a_j}$.
This makes the analytic morphism $\Stm$ from $\Mod(F)=\Spec(\mathbb{C}[a_1,a_2,a_3]/(a_1+a_2+a_3))$ to
$\Stokesmoduli(F)=\Spec(\mathbb{C}[x_{0,1},x_{2,1}])$ explicit.

{\it Special cases:}
\begin{enumerate}\itemsep=0pt
\item[(1)] $a_1=-1/3$, $a_2=0$, $a_3=1/3$ yields $x_{0,1}=x_{2,1}=0$ and all Stokes matrices are trivial.

\item[(2)] $a_1=a_2=a_3=0$ yields $x_{0,1}={3\choose 1}$, $x_{2,1}=-{3\choose 2}$.
The equation is $\delta^3-z$ and all Stokes entries are $\pm $ binomial coef\/f\/icients.
\end{enumerate}

{\it For general~$n$}, $\Mod(F)$ is represented by the universal family
\begin{gather*}
\delta +\left(
\begin{matrix} a_0&0 & \ddots& \ddots &z
\\
1 &a_1 &0 &\ddots & 0
\\
\vdots &\vdots &\ddots &\ddots & \ddots
\\
\cdot&\dots & 1&a_{n-2} &0
\\
0& \dots &\dots & 1& a_{n-1}
\end{matrix}
\right)
\qquad
\text{with}
\qquad
\sum a_j=0.
\end{gather*}
The case $a_i=\frac{i}{n}-\frac{n-1}{2n}$ for $i=0,\dots,n-1$ corresponds to the case where all Stokes matrices are trivial.

The case all $a_j=0$ corresponds to $\delta^n-z$.
For this equation, as in the case $n=3$, all the Stokes entries are $\pm$ binomial coef\/f\/icients.
This is an explicit form of part of a~conjecture of Dubrovin (see Section~\ref{section3}).

{\bf Conclusions.} For the totally ramif\/ied cases~$F$ that we consider, the Stokes matrices have explicit formulas in
exponentials of algebraic expressions in the entries of the matrix dif\/ferential operator.
This shows in particular, that there is no universal family parametrized by $\Stokesmoduli(F)$.
Further $\Stm: \Mod(F)\rightarrow \Stokesmoduli(F)$ is surjective and the f\/ibers correspond to choices of the logarithm of
the topological monodromy.

Also for a~slightly more general case than the totally ramif\/ied case the entries the Stokes matrices are determined~by
the topological monodromy.

\begin{proposition}\label{proposition2}
Let~$M$ be an irreducible differential module over the field $K=\mathbb{C}(\{z^{-1}\})$.
Suppose that the tuple $(V,\dots)$, associated by Theorem~{\rm \ref{theorem2}} to~$M$, has the properties:
\begin{enumerate}\itemsep=0pt
\item[$(a)$]
$V=V_0\oplus W$ where $V_0$ has eigenvalue~$0$, $\dim V_0=m$, and~$W$ is totally ramified, $\dim W=n$.
\item[$(b)$]
The restriction of the formal monodromy~$\gamma$ to $V_0$ has~$m$ distinct eigenvalues.
\end{enumerate}

Then, after normalization, the monodromy identity implies that the topological monodromy at $z=\infty$ determines all
Stokes matrices.
\end{proposition}

\begin{proof}
The Malgrange irregularity $\irr(M)$ is $2m+n-1$.
Let $f_1,\dots,f_m$ denote a~basis of eigenvalues of~$\gamma$ on $V_0$.
Further we may suppose that $W=\oplus_{i=0}^{n-1}V_{\zeta^iz^{1/n}}$, where $\zeta =e^{2\pi i/n}$ and each $V_{\zeta^iz^{1/n}}$ has dimension 1.
Let $e_i$ be a~basis vector for $V_{\zeta^iz^{1/n}}$ such that $\gamma e_i=e_{i+1}$ for $i=0,\dots,n-2$ and $\gamma
e_{n-1}=\lambda e_0$ for some $\lambda \in \mathbb{C}^*$.

The Stokes maps $\operatorname{St}_d$ for directions $d\in [0,1)$ are maps $\ell_i:V_0\rightarrow V_{\zeta^iz^{1/n}}$,
$\tilde{\ell}_i:V_{\zeta^iz^{1/n}}\rightarrow V_0$ and $m_{i,j}:V_{\zeta^iz^{1/n}}\rightarrow v_{\zeta^jz^{1/n}}$ for $i\neq j$.
Write $\ell_i(f_j)=x_{i,j}e_i$, $\tilde{\ell}_i(e_i)=\sum\limits_jy_{j,i}f_j$ and $m_{i,j}(e_i)=z_{i,j}e_j$.

An element $c=(c_1,\dots,c_m)\in (\mathbb{C}^*)^m$ acts on the basis $f_1,\dots,f_m$ of $V_0$ by $f_i\mapsto c_if_i$ for all~$i$.
Then~$c$ acts on the Stokes entries by $x_{i,j}\mapsto c_jx_{i,j}$, $y_{j,i}\mapsto c_j^{-1}y_{j,i}$ and $z_{i,j}\mapsto z_{i,j}$.

Suppose that a~$c=(c_1,\dots,c_m)\neq 1$ acts trivially on the above Stokes entries.
Then $c_j\neq 1$ implies $x_{i,j}=y_{j,i}=0$.
Let $\tilde{V}_0$ be the subspace of $V_0$ generated by the $f_j$ with $c_j=1$.
Then $(\tilde{V}_0\oplus W,\dots)$ is a~subobject of $(V,\dots)$ and, as a~consequence,~$M$ is reducible.

Thus the action of $(\mathbb{C}^*)^m$ is faithful and one can produce a~basis $f_1,\dots,f_m$ such that~$m$ of the
Stokes entries are~1.
This leaves $m+n-1$ unknown Stokes entries.
The characteristic polynomial of the product $\gamma \prod\limits_{d\in [0,1)} \operatorname{St}_d$ can be explicitly computed
(compare~\cite{CM-vdP}) and the entries of this polynomial are independent inhomogeneous linear expressions in the
$z_{i,j}$ and the products $x_{i,j}y_{j,i}$.
Thus the monodromy identity and normalizing some of the $x_{i,j}$, $y_{j,i}$ to 1 produces all $\operatorname{St}_d$ for $d\in[0,1)$.
For general~$d$, the Stokes matrix $\operatorname{St}_d$ is derived from the above, using the identity $\gamma^{-1} \operatorname{St}_d\gamma=\operatorname{St}_{d+1}$.
\end{proof}

\section{Fano varieties and quantum dif\/ferential equations}\label{section3}

This part of the paper reviews work by J.A.~Cruz Morales and the author~\cite{CM-vdP}.
A~(complex) Fano variety~$F$ is a~non-singular, connected projective variety of dimension~$d$ over $\mathbb{C}$, whose
anticanonical bundle $(\Lambda^d\Omega)^*$ is ample.
There are rather few Fano varieties.

{\it Examples}: For dimension $1$ only $F=\mathbb{P}^1$; for dimension 2: the Fano's are del Pezzo surfaces
and $\cong \mathbb{P}^1\times \mathbb{P}^1$ or $\cong$ to $\mathbb{P}^2$ blown up in at most 8 points in general
position; for dimension~3 the list of Iskovshih--Mori--Mukai contains 105 deformation classes.
According to Wikipedia, the classif\/ication of Fano varieties of higher dimension is a~project called ``the periodic
table of mathematical shapes''.

\subsection{Quantum cohomology and quantum dif\/ferential equations}\label{section3.1}

We borrow from the informal introduction to the subject from M.A.~Guest's book~\cite{Gue}.
Let~$F$ be a~Fano variety.
On the vector space $H^*(F):=\oplus_{i=0}^dH^{2i}(F,\mathbb{C})$ there is the usual cup product $\circ $ (say obtained
by the wedge product of dif\/ferential forms).
Quantum cohomology introduces a~deformation $\circ_t$ of the cup product $\circ $ on $H^*(F)$ for $t\in H^2(F,\mathbb{C})$.

With respect to a~basis $b_0,\dots,b_s$ of $H^*(F,\mathbb{Z}):=\oplus_{i=0}^dH^{2i}(F,\mathbb{Z})$, the quantum products
$b_i\circ_t $ have matrices which are computable in terms of the geometry of~$F$.
Let $b_1,\dots,b_r$ be a~basis of $H^2(F,\mathbb{Z})$.

{\it The quantum differential equation of a~Fano variety~$F$} is a~system of (partial) linear dif\/fe\-ren\-tial operators
$\partial_i-b_i\circ_t$, $i=1,\dots, r$ acting on the space of the holomorphic maps $H^2(F,\mathbb{C})=\mathbb{C}b_1+\dots
+\mathbb{C}b_r\rightarrow H^*(F,\mathbb{C})= \mathbb{C}b_0+\dots +\mathbb{C}b_s$.

For the case $r=1$ that interests us, the quantum dif\/ferential equation reads $z\frac{d}{dz}\psi=C\psi$, where~$\psi$ is
a~vector of lenght $s+1$ and~$C$ is the matrix of quantum multiplication $b_1\circ_t$.
The entries of the $(s+1)\times (s+1)$ matrix~$C$ are polynomials in~$z$ with integer coef\/f\/icients.
Clearly $z=0$ is a~regular singular point and $z=\infty$ is irregular singular.
By taking a~cyclic vector one obtains a~scalar dif\/ferential equation of order $s+1$.

\subsection{Examples of quantum dif\/ferential equations}\label{section3.2}

Below we present examples of quantum dif\/ferential equations.

$\delta^n-z$ with $\delta=z\frac{d}{dz}$ for $\mathbb{P}^{n-1}$.

$\delta^{n+m-1}-m^mz(\delta +\frac{m-1}{m})(\delta +\frac{m-2}{m})\dots (\delta +\frac{1}{m})$ with $n\geq 1$, $m>1$
for a~non-singular hypersurface of degree~$m$ in $\mathbb{P}^{n+m-1}$~\cite[\S~3.2, Example~3.6, p.~43]{Gue}.

$\prod\limits_{j=0}^n\delta \big(\delta -\frac{1}{w_j}\big)\cdots \big(\delta -\frac{w_j-1}{w_j}\big) -z$ for the weighted projective
space $\mathbb{P}(w_0,\dots,w_n)$ (an orbifold).

$\delta^3-az\delta^2-((b-a)z^2+bz)\delta +2az^2-cz^3$ with $a,b,c\in \mathbb{Z}$ for del Pezzo surfaces.

$\delta^4-11z\delta^2-11z\delta -3z-z^2$ for $V_5$,
for a~linear section of the Grasmannian $G(2,5)$ in the Pl\"ucker embedding.

$\delta^4-(94z^2+6z)\delta^2-(484z^3+188z^2+2z)\delta -(695z^4+632z^3+98z^2)$ for the 3-fold $V_{22}$.

B.~Dubrovin is one of the founders of quantum cohomology.
One of his {\it conjectures}~\cite{Du98, Dub99} states that the Gram matrix $(G_{i,j})$ of a~(good) Fano
variety coincides with the ``{\it Stokes matrix}'' of the quantum dif\/ferential equation of~$F$ (up to a~certain
equivalence of matrices).

Here $G_{i,j}=\sum\limits_k (-1)^k \dim \operatorname{Ext}^k(\mathcal{E}_i,\mathcal{E}_j)$, where $\{\mathcal{E}_i\}$ is an exceptional
collection of coherent sheaves on~$F$ generating the derived category $D^b\operatorname{coh}(F)$.
Further ``Stokes matrix'' is in fact a~connection matrix and, in our terminology, equal to the product
$\prod\limits_{d\in [0,1/2)}\operatorname{St}_d$ (in counter clock order).

For $F=\mathbb{P}^n$, this has been verif\/ied by D.~Guzzetti~\cite{Guz}.
There are recent papers~\cite{Guestetal,GuestSakai,Iri,MT,Tan,T-U,Ue1,Ue2} which handle more cases, e.g., Calabi--Yau complete intersections
in weighted projective spaces, Grassmanians, cubic surfaces, toric orbifolds.
They use the ``Fourier--Laplace type transformation'' to a~regular singular dif\/ferential equation, mentioned in Section~\ref{section2.1}.

The contribution of~\cite{CM-vdP} is proving Dubrovin's conjecture by computing all $\operatorname{St}_d$, using only the formal
classif\/ication and {\it the monodromy identity}, for the cases $\mathbb{P}^n$, non singular hypersurfaces of degree
$m\leq n$ in $\mathbb{P}^n$, and for weighted projective spaces $\mathbb{P}(w_0,\dots,w_n)$.
The equations and the method are closely related to Sections~\ref{section1.5} and \ref{section2.2}.

\section{Riemann--Hilbert approach to Painlev\'e equations}\label{section4}

This classical method, related to isomonodromy, was revived and ref\/ined by M.~Jimbo, T.~Miwa and K.~Ueno~\cite{JMU, JM}.
The literature on the subject is nowadays impressive.

{\sloppy Certain details of the Riemann--Hilbert approach are worked out  in collaboration with \mbox{M.-H.~Saito}~\cite{vdP-Sa}.
In collaboration with J.~Top, ref\/ined calculations of Okamoto--Painlev\'e spaces and B\"acklund transformations were
presented for ${\rm P_{I}}$--${\rm P_{IV}}$ (2009--2014), see~\cite{vdP-T} and its references.
We give here a~{\it rough sketch} of the ideas and especially of the part where Stokes matrices enter the picture.
The starting point is a~`family'
$\bf S$ of dif\/ferential modules~$M$ over~$\mathbb{C}(z)$ with prescribed singularities
at f\/ixed points of~$\mathbb{P}^1$.
A~priori $\bf S$ is just a~set.
One side of the Riemann--Hilbert approach is to construct a~moduli space~$\mathcal{M}$ over~$\mathbb{C}$ such that~$\bf S$
has a~natural identif\/ication with $\mathcal{M}(\mathbb{C})$.
This part does not involve Stokes maps.

}

The other side of the Riemann--Hilbert approach is a~``monodromy space'' $\mathcal{R}$ built out of monodromy, Stokes matrices and `links'.
The space $\mathcal{R}$ is determined by the prescribed type and position of the singularities of $\bf S$.

{\it An example for $\mathcal{R}$}.
The set $\bf S$ which gives rise to ${\rm P_{III}}(D7)$ consists of the dif\/ferential modu\-les~$M$ over $\mathbb{C}(z)$
which have only $0$ and~$\infty$ as singular points.
The point $0$ has Katz invariant $1/2$ and the point~$\infty$ has Katz invariant $1$ (see Section~\ref{section1.2} for the def\/inition of
the Katz invariant).
The monodromy space $\mathcal{R}$ consists of the analytic classif\/ication $(V(0),\dots)$ of~$M$ at $z=0$ and
$(V(\infty),\dots)$ at $z=\infty$ and a~connection matrix between these data, the link $L:V(0)\rightarrow V(\infty)$,
which describes the relation between the solutions around $z=0$ and the solutions around $z=\infty$.

The solution space $V(0)$ at $z=0$ is given a~basis $e_1$, $e_2$ for which the formal monodromy, the Stokes matrix and
topological monodromy $\operatorname{top}_0$ are
\begin{gather*}
\begin{pmatrix} 0 & -1\\ 1 &  0\end{pmatrix},\qquad \begin{pmatrix} 1 & 0\\ e & 1\end{pmatrix}, \qquad \begin{pmatrix}
-e & -1 \\ 1 &  0\end{pmatrix}.
\end{gather*}
The solution space $V(\infty)$ at $z=\infty$ is given a~basis $f_1$, $f_2$ for which the formal monodromy, the Stokes maps
and the topological monodromy $\operatorname{top}_\infty$ are
\begin{gather*}
\begin{pmatrix} \alpha & 0 \\ 0 & \alpha^{-1} \end{pmatrix},
\qquad \begin{pmatrix} 1 &  0\\  c_1 & 1\end{pmatrix}, \qquad \begin{pmatrix} 1 & c_2\\ 0 & 1\end{pmatrix},\qquad
\begin{pmatrix}\alpha & \alpha c_2
\\
\alpha^{-1}c_1&\alpha^{-1}(1+c_1c_2)
\end{pmatrix}.
\end{gather*}
One may assume that $L:=\begin{pmatrix} \ell_1 & \ell_2\\ \ell_3 & \ell_4\end{pmatrix}$ has determinant~1.
There is a~relation $\operatorname{top}_0\cdot L^{-1}\cdot \operatorname{top}_\infty \cdot L=1$.
This yields a~set of variables and relations and thus an af\/f\/ine variety.
The above bases~$e_1$,~$e_2$ of~$V(0)$ and $f_1$, $f_2$ of~$V(\infty)$ are not unique.
Indeed, the ambiguity in these basis is given by the transformation $e_1,e_2\mapsto \lambda_0e_1, \lambda_0e_2$
and $f_1,f_2\mapsto \lambda_1f_1, \lambda_2f_2$
with $(\lambda_0,\lambda_1,\lambda_2)\in (\mathbb{C}^*)^3$.
By dividing this af\/f\/ine space by the action of $(\mathbb{C}^*)^3$ one obtains $\mathcal{R}$.
The f\/inal result is that {\it $\mathcal{R}$ is an affine cubic surface, given by variables $\ell_{13}$, $\ell_{23}$,
$\alpha$ and relation $ \ell_{13}\ell_{23}e+\ell_{13}^2+\ell_{23}^2+\alpha \ell_{13}+\ell_{23} =0$.}

Consider the map ${\bf S}\rightarrow \mathcal{R}$ which associates to each module $M\in {\bf S}$ its monodromy data in
$\mathcal{R}$.
The f\/ibers of this map are parametrized by some $T\cong \mathbb{C}^*$ and there results a~bijection ${\bf S}\rightarrow
\mathcal{R}\times T$.
The set $\bf S$ has a~priori no structure of an algebraic variety.
A~moduli space $\mathcal{M}$ over $\mathbb{C}$, whose set of closed points consists of certain connections of rank two
on the projective line, is constructed such that $\bf S$ coincides with $\mathcal{M}(\mathbb{C})$.
This def\/ines the analytic Riemann--Hilbert morphism $\operatorname{RH}\colon \mathcal{M}\rightarrow \mathcal{R}$.
The f\/ibers of $\operatorname{RH}$ are the isomonodromic families.
There results an extended Riemann--Hilbert {\it isomorphism} $\operatorname{RH}^+\colon \mathcal{M}\rightarrow \mathcal{R}\times T$.
From the isomor\-phism~$\operatorname{RH}^+$ the Painlev\'e property for the corresponding Painlev\'e equation follows and the moduli
space $\mathcal{M}$ is identif\/ied with an Okamoto--Painlev\'e space.
Special properties of solutions of the Painlev\'e equations, such as special solutions, B\"acklund transformations etc.,
are derived from the extended Riemann--Hilbert isomorphism.

The above sketch needs subtle ref\/inements.
One has, depending on the Painlev\'e equation and its parameters, to add level structure, to forget points, to
desingularize $\mathcal{R}$ and $\mathcal{M}$, to replace spaces by their universal covering etc., in order to obtain
a~correct extended Riemann--Hilbert isomorphism.
For the remarkable fact that for each Painlev\'e equation the moduli space for the monodromy $\mathcal{R}$ is an af\/f\/ine
cubic surface with three lines at inf\/inity, there is not yet an explanation.

\subsection*{Acknowlegdements}
The author likes to thank the referees for their work and the very useful comments which led to many improvements.

\pdfbookmark[1]{References}{ref}
\LastPageEnding

\end{document}